\documentclass[12pt]{amsart}
\usepackage{amssymb,amsmath,amscd}
\newtheorem{thm}{Theorem}[section]
\newtheorem{corollary}[thm]{Corollary}

\newtheorem{lemma}[thm]{Lemma}
\newtheorem{fact}[thm]{Fact}

\theoremstyle{definition}

\newtheorem{example}[thm]{Example}

\newtheorem{remark}[thm]{Remark}

\newcommand{\ra}{\rightarrow}

\newcommand{\Pthree}{\mathbb{P}^3}
\newcommand{\ptwo}{\mathbb{P}^2}

\begin{document}

\title{Holomorphic injectivity and the Hopf map}

\author{Scott Nollet}
\address{Department of Mathematics,
Texas Christian University,
Fort Worth, TX 76129}
\email{s.nollet@tcu.edu}

\author{Frederico Xavier}
\address{Department of Mathematics,
University of Notre Dame, Notre Dame, IN 46556}
\email{fxavier@nd.edu}

\date{}
\dedicatory{Warmly dedicated to Joana}
\thanks{Work of the second author was partially supported by NSF 
grant DMS02-03637}

\begin{abstract}
    We give sharp conditions on a local biholomorphism 
    $F:X \to \mathbb C^{n}$ which ensure global injectivity.
    For $n \geq 2$, such a map is injective if for each complex line 
    $l \subset \mathbb C^{n}$, the pre-image $F^{-1}(l)$ 
    embeds holomorphically as a connected domain into 
    $\mathbb C \mathbb P^{1}$, the embedding being 
    unique up to M\"obius transformation. 
    In particular, $F$ is injective if the pre-image of every complex 
    line is connected and conformal to $\mathbb C$. 
    The proof uses the topological fact that the natural map $\mathbb R 
    \mathbb P^{2n-1} \to \mathbb C \mathbb P^{n-1}$ associated to the 
    Hopf map admits no continuous sections and the classical Bieberbach-Gronwall 
    estimates from complex analysis.
\end{abstract}

\maketitle

\section{Introduction}

The study of univalence (injectivity) of holomorphic functions of one variable is a 
classical topic in complex analysis. 
One of the high points in the subject was the celebrated solution of the Bieberbach 
Conjecture by de Branges in 1984 \cite{deB,gong,hayman}. 
As part of the effort to understand the conjecture, several authors introduced 
various univalence criteria for locally univalent holomorphic functions defined on 
the open unit disc in $\mathbb C$ \cite{pommerenke}. 
In contrast, fewer injectivity criteria are known for local biholomorphisms 
in higher dimensions \cite[chapter V]{gong}.
In the simplest case when $F: \mathbb C^n \to \mathbb C^n$ is a {\it polynomial} 
local biholomorphism, it is not even known whether $F$ is a priori injective 
(hence bijective \cite[I, Theorem. 2.1]{BCW}) without further hypotheses: 
this is the jacobian conjecture (see \cite{BCW,E} for general references), 
which remains open after more than 60 years.

In this paper we draw a connection between global injectivity of a local 
biholomorphism $F$ and connectedness of certain pre-images of $F$. 
The na\"ive observation that $F$ is injective if and only if the pre-image of 
each point is connected leads one naturally to ask whether there might be 
a similar criterion based on connectedness of pre-images of {\it positive} dimensional submanifolds. 
For example, it's easy to see that a local diffeomorphism 
$F:X \to \mathbb R^{n}$ 
is injective if for each real line $l$, the pre-image $F^{-1}(l)$ is connected 
\cite[Example 3.3]{NX}.
For a local biholomorphism $F:X \to \mathbb C^{n}$ and {\it complex} lines 
$l \subset \mathbb C^{n}$, it turns out that  connectedness of the pre-images is 
not enough (Example \ref {highergenus}), one also needs information about 
the conformal type of the pre-images $F^{-1}(l)$ 
(Corollaries\ref {global} and \ref {algebraic1}). 

To state our injectivity criterion, we need the concept of a 
{\it rigid domain} of $\mathbb C \mathbb P^1$: 
for us, these are the connected open sets 
$i:U \hookrightarrow \mathbb C \mathbb P^1$ such that any holomorphic 
embedding $f:U \hookrightarrow \mathbb C \mathbb P^1$ differs from 
$i$ by an automorphism $M$ of $\mathbb C \mathbb P^1$, 
that is $f = M \circ i$.
A rigid domain $U \subset \mathbb C \mathbb P^1$ is necessarily dense
(otherwise apply an automorphism of $\mathbb C \mathbb P^{1}$ that 
takes $U$ into the unit disc $D \subset \mathbb C$. 
By the Riemann mapping theorem, there are many holomorphic embeddings 
$F:D \hookrightarrow \mathbb C \subset \mathbb C \mathbb P^1$ which are 
not the restriction of a M\"obius transformation. 
By unicity of continuation, the restriction $F|_{U}$ is not either). 
Typical examples include $\mathbb C \subset \mathbb C \mathbb P^{1}$, 
the complement of $\mathbb C$ by finitely many points 
(a connected rational curve, in the sense of algebraic geometry) 
or, more generally, the complement in the Riemann sphere of a removable closed set 
(a closed subset $E \subset \mathbb C\mathbb P^1$ is {\it removable} if for 
each open set $U \subset \mathbb C \mathbb P^1$, the bounded holomorphic 
functions $f:U-E \to \mathbb C$ extend holomorphically to $U$).

\begin{thm}\label{main}
     Let $X$ be a connected complex manifold of dimension $n \geq 2$, 
     $F:X \to \mathbb C^{n}$ a local biholomorphism. 
     Fix $q \in F(X)$ and suppose that $F^{-1}(l)$ is conformal to a 
     rigid domain $D_{l} \subset \mathbb C \mathbb P^{1}$ for every complex 
     line $l$ passing through $q$. 
     Then $q$ is assumed exactly once by $F$.
\end{thm}

The statement is sharp in all respects, starting with the obvious condition 
$n\geq 2$. 
Its conclusion may fail if even a pre-image of a single line is disconnected 
(Example \ref{single-analytic}) or if the pre-images are conformal to punctured 
compact Riemann surfaces of {\it positive} genera (Example 
\ref{highergenus}). 
The rigidity condition also cannot be removed (Example \ref{rigid}). 
The proof is a blend of analytic and topological ideas: 
if $F^{-1}(q)$ contains {\it two} points, we can 
construct a {\it continuous} section to the natural map
$\pi: \mathbb R \mathbb P^{2n-1} \to \mathbb C \mathbb P^{n-1}$ 
associated to the Hopf map, an impossibility. The continuity of this 
section follows from the classical Bieberbach-Gronwall estimates for univalent 
functions on the unit disc. 
Globalizing Theorem \ref{main} yields:

\begin{corollary}\label{global}
Let $F:X \to \mathbb C^{n}$be a local biholomorphism with $n \geq 2$.
If each non-empty pre-image of a complex line is connected and 
conformal to $\mathbb C$, then $F$ is injective.
\end{corollary}
 %       A local biholomorphism $F:X \rightarrow \mathbb C^{n}$ 
%        satisyfing $n \geq 2$ is 
%        injective if the pre-image of every complex line is 
%	connected and conformal to an open rigid domain in 
%	$\mathbb C \mathbb P^{1}$.

Applying Serre's GAGA principle \cite{serre} gives an algebro-geometric
variant:

\begin{corollary}\label{algebraic1}
         Let $F:X \ra \mathbb A^{n}_{\mathbb C}$ be an \'etale morphism
         of schemes with $n \geq 2$. Then the following are equivalent.
\begin{enumerate}
         \item[(1)] $F$ is injective.
         \item[(2)] For each line $l \subset \mathbb C^{n}$ meeting $F(X)$, 
	 the pre-image $F^{-1}(l)$ is connected and rational.
\end{enumerate}
\end{corollary}
\noindent The forward direction is clear because the pre-image of a 
line $l$ is identified with a Zariski open subset of $l$; the 
converse is immediate from Theorem \ref{main} above. 
Example \ref{highergenus} shows that rationality is needed in condition (2) above. 
 
Finally, applying generic smoothness gives a projective version: 

\begin{corollary}\label{algebraic2}
         Let $F:Z \ra \mathbb P^{n}_{\mathbb C}$ be a generically finite morphism 
         of schemes with $n \geq 2$ and $Z$ smooth. Then the following are equivalent.
\begin{enumerate}
         \item[(1)] $F$ is birational.
         \item[(2)] There is an open set $U\subset \mathbb P^{n}$ 
         such that for each line $l \subset \mathbb P^{n}$, the 
         pre-image $F^{-1}(l \cap U)$ is either empty or irreducible and 
         rational.
\end{enumerate}
\end{corollary}

\begin{proof}
    Since $\mathbb C$ has characteristic zero, 
    generic smoothness \cite[III, Cor. 10.7]{AG} gives an 
    open set $V \subset \mathbb P^{n}$ for which the restricted map 
    $F:F^{-1}(V) \to V$ is smooth. Generic finiteness of $F$ implies 
    that (a) the set $V$ is non-empty (hence dense)  and (b) the 
    relative dimension is zero so that $F:F^{-1}(V) \to V$ is \'etale. 
    Further restriction to a standard open affine subset 
    $\mathbb A^{n} \cong U_{i} \subset \mathbb P^{n}$ 
    gives an \'etale morphism $F^{-1}(V \cap \mathbb A^{n}) \to V \cap \mathbb 
    A^{n} \hookrightarrow \mathbb A^{n}$. \\
   $\Rightarrow:$ If $F$ is birational, then $\deg F = 1$ and we may take 
   $U=V \cap \mathbb A^{n}$. \\
   $\Leftarrow:$ The second condition continues to hold if we 
   replace $U$ with $U \cap V \cap \mathbb A^{n}$. Composing with 
   the inclusion into $\mathbb A^{n}$, Corollary \ref{algebraic1} 
   implies that $F|_{F^{-1}(U \cap V \cap \mathbb A^{n})}$ is injective, 
   hence $F$ is birational.
\end{proof}
 
In section two we prove Theorem \ref{main}, modulo a
delicate continuity proof. We also give examples which show that the 
hypotheses cannot easily be removed.  
Section three is devoted to the continuity proof.
We direct the interested reader to \cite{meisters,NX,parth,SX,X1} and 
the references therein for further results on global injectivity in the
differentiable context.
\vskip10pt
\noindent
{\bf Acknowledgements.} In an earlier version of this paper, the conclusion of 
Theorem \ref{main} was that the point $q$ can be assumed at most 
\it twice \rm by $F$ and there was an example intended to show that $q$ could 
have two points in its pre-image. 
J\'anos Koll\'ar generously pointed out how the example failed and 
suggested in broad terms how our method could be strengthened to yield the 
optimal result. 
He also indicated that Corollary \ref{algebraic1} can be approached using ideas 
from Mori theory \cite{kollar,mori}.

We would also like to thank Francis Connolly, Richard Hind, and Sean Keel for 
useful conversations. 

\section{Main ideas and examples}\label{markers}

In this section we prove Theorem \ref{main} except for a lengthy 
technical matter which is delayed to section three. 
We also give examples to show that the hypotheses cannot be weakened.

Recall that the Hopf map $S^{2n-1} \to \mathbb C \mathbb P^{n-1}$ sends a 
unit vector $u \in \mathbb C^{n}$ to the complex one-dimensional subspace
containing it. 
Clearly this induces a map 
$\pi: \mathbb R \mathbb P^{2n-1} \to \mathbb C \mathbb P^{n-1}$ 
between real and complex projective spaces. 

\begin{fact}\label{nosections}
    Neither $\pi$ nor the Hopf map admits a continuous section for $n \geq 2$. 
\end{fact}

For instance, the composite map in cohomology
$$H^2 (\mathbb C \mathbb P^{n-1})\rightarrow H^2(\mathbb R \mathbb P^{2n-1})
\rightarrow H^2 (\mathbb C \mathbb P^{n-1})$$
induced by a continuous section of $\pi$ must be the identity, but 
this is impossible since $H^2 (\mathbb C \mathbb P^{n-1})\neq 0 $ and 
$H^2(\mathbb R \mathbb P^{2n-1})=0$.

This topological fact can sometimes be used to prove injectivity of 
local biholomorphisms. The following example illustrates the idea. 

\begin{example}\label{geodesics}
        Let $F:  \mathbb C^{n} \rightarrow \mathbb C^n$ be a local
	bihilomorphism with $n \geq 2$ such that the pre-image of every
	complex line is connected and simply connected.
	Then $F$ is injective.
\end{example}

\begin{proof} If $F$ is not injective, we may suppose that $F(p)=F(q)=0$
with $p \neq q$, hence $F^{-1}(l)$ contains both $p$ and $q$ for all
one-dimensional complex subspaces $l$.
It follows from the inverse function theorem that the complex curve
$F^{-1}(l)$ is properly embedded in $\mathbb C^{n}$ (whether $F$ is a
proper map or not), hence with respect to the induced Riemannian metric it
is a complete simply connected real surface of non-positive curvature
\cite[p. 79]{gh}.

It follows from Hadamard's theorem \cite[Ch. 7, Theorem 3.1]{dC}
that any two points in $F^{-1}(l)$ can be joined by a {\it unique} geodesic of 
$F^{-1}(l)$.
Given $l \in \mathbb P^{n-1}$, let $w(l)$ denote the
initial vector of the (unique) unit-speed geodesic along $F^{-1}(l)$
joining $p$ to $q$ and set $v(l)=dF(0) w(l) \in T_{l,0}$.
Notice that all geodesic segments are contained in a fixed compact
subset of $\mathbb C^{n}$.
The map $v$ is continuous because geodesics converge to geodesics in the
$C^2$ topology (which, after passing to subsequences, is a consequence of
uniform $C^3$ boundedness) and from the uniqueness of the geodesic
along $F^{-1}(l)$.
Since $v$ is nonvanishing, it is clear that $\frac{v(l)}{|v(l)|}$ 
defines a continuous section to the Hopf map, a contradiction in view 
of \ref{nosections}.
\end{proof}

\begin{remark}
    Notice in the previous example that if the pre-images of complex 
    lines are not simply connected, then the geodesics along them 
    from $p$ to $q$ may not be unique, so the construction above does 
    not work. On the other hand, because the curvature is 
    non-positive, there is a unique geodesic for each 
    homotopy class of paths from $p$ to $q$, so there is at least a 
    local finite multisection in this case: potentially pieces of these 
    could be glued to construct a global continuous section. We 
    conjecture that for local biholomorphisms 
    $F:\mathbb C^{n} \to \mathbb C^{n}$ that connectedness of the pre-images 
    alone is enough to ensure injectivity of $F$. 
\end{remark}

The proof of our main theorem is similar in spirit to Example \ref{geodesics}, except 
that we use tangent vectors to rational curves to produce a section to the map $\pi$. 
The construction of the map is given below, the proof of continuity 
being delayed to section three. 

\vskip10pt \noindent{\bf Proof of Theorem \ref{main}}
\vskip5pt       
Taking $q=0$, suppose that $F^{-1}(0)$ contains at least two distinct points 
$z_{1} \neq z_{2}$. 
Fix the two points $w_{1}=0$ and $w_{2}=1$ in $\mathbb C \mathbb P^1$ 
along with two nonzero real tangent vectors 
$0 \neq \alpha_i \in T_{\mathbb C \mathbb P^{1},w_{i}}$.
Let $l$ be a complex line through $0$. 
We claim that 

\vskip10pt

\noindent (i) There are exactly two holomorphic embeddings 
$T_{l}=T_l^{1}, T_{l}^{2}:F^{-1}(l)\to \mathbb C \mathbb P^1$ such that 
$T_{l}(z_{i})=w_{i}$ ($i=1,2$) and 
\begin{equation}\label{uno}
    dF(z_1)(dT_{l}(z_1))^{-1}\alpha_1=dF(z_2)(dT_{l}(z_2))^{-1}\alpha_2.
\end{equation}  

\vskip10pt

\noindent (ii) Letting $v_{1}$ (resp. $v_{2}$) be the vector in 
equation (\ref{uno}) for $T_{l}^{1}$ (resp. $T_{l}^{2}$), we have $v_{1}=-v_{2}$. 

\vskip10pt

To see this, fix a holomorphic embedding 
$S_l:F^{-1}(l) \hookrightarrow \mathbb C \mathbb P^{1}$ 
that takes $z_i$ to $w_i$. By rigidity, any such embedding 
$F^{-1}(l)\to \mathbb C \mathbb P^1$ is of the form $T_l=U \circ S_l$ for an 
automorphism $U: \mathbb C \mathbb P^1 \to \mathbb C \mathbb P^1$ 
which fixes $w_{1}=0$ and $w_{2}=1$. Any such $U$ has the form 
$$U_a(z)= \frac{z}{az+1-a}$$ for some $a\neq 1$.
Setting $b=1-a$, and using $U_a'(0)= b^{-1}$, 
$U_a'(1)=b$ we substitute $(dT_l(z_1))^{-1}\alpha_1= bdS_l^{-1}(0) \alpha_1$
and $(dT_l(z_2))^{-1}\alpha_2= b^{-1} dS_l^{-1}(1) \alpha_2$ into 
Equation (\ref{uno}) to obtain
$$bdF(z_1) dS_l^{-1}(0) \alpha_1 = b^{-1}dF(z_2) dS_l^{-1}(1) \alpha_2$$
in the tangent space $T_{l,0} \subset T_{\mathbb C^{n},0}$, which we 
identify with the subset $\mathbb C \cong l \subset \mathbb C^{n}$. 
Thus $b^{2} \in \mathbb C - \{0\}$ and we find two values for $b$, one being 
the negative of the other. 
This proves (i) and (ii) above. 

In particular, both choices of $T_l$ yield the same {\it real} line 
$v(l) \subset l \subset \mathbb C^n$, hence the map 
$v:\mathbb C \mathbb P^{n-1} \to \mathbb R \mathbb P^{2n-1}$ is a 
set-theoretic section of the map $\pi$ above. 
Using the Bieberbach-Gronwall estimate for univalent functions on 
the unit disc \cite[Theorem 1.3]{hayman}, we will show in section three 
that $v$ is a {\it continuous} section: this contradiction finishes the proof. 
\hfill \qed

\vskip10pt

\begin{example}\label{single-analytic}
The conclusion of Theorem \ref{main} may fail if the pre-image
of even a single line is disconnected.
Let $f:\mathbb C \to \mathbb C$  be a holomorphic function with nowhere zero
derivative and which assumes the value $0$ more than once.
Let $F$ be the local biholomorphism of
$\mathbb C^{2}$ defined by $F(z,w)= (f(z),w)$.
If $l \subset \mathbb C^2$ is the complex line $aw_1+bw_2=0$, then
$F^{-1}(l)$ is conformal to $\mathbb C$ if $b\neq 0$, which is a rigid domain. 
For $b=0$, however, the pre-image $F^{-1}(l)$ is homeomorphic to the 
disconnected space $f^{-1}(0) \times \mathbb C$.
\end{example}

\begin{example}\label{highergenus}
    This example shows that the rationality condition in part (2) of Corollary 
    \ref{algebraic1} cannot be removed, hence the condition that 
    $F^{-1}(l)$ be conformal to a rigid domain of $\mathbb C \mathbb P^1$, 
    rather than a compact Riemann surface of positive genus, cannot be removed 
    from Theorem \ref{main} or Corollary \ref{global}.
    Let $Y \subset \Pthree$ be a general smooth surface of degree $d > 2$.
    Then $Y$ is neither ruled nor a Steiner surface, hence has no
    two-dimensional families of reducible plane sections 
    \cite[Lemma II.2.4]{al}.
    It follows that for a general point $p \in \Pthree$, only a {\it finite} 
    number of planes $H_{i}$ containing $p$ intersect $Y$ in a reducible 
    curve. 
    Let $F:Y \to \ptwo$ be the projection from such a point $p$ onto a 
    plane $\ptwo \subset \Pthree$.
    If $q \in \ptwo$ and $q \not \in \bigcup_{i} H_{i}$, then {\it every} line
    $l$ through $q$ has irreducible pre-image $F^{-1}(l) = Y \cap H$ 
    ($H$ spanned by $l$ and $p$) of arithmetic 
    genus $p_{a}(Y \cap H)=\frac{1}{2}(d-1)(d-2)$ \cite[I, Ex. 7.2(b)]{AG}. 
    Moreover, $Y \cap H$ is smooth and irreducible for general 
    $l$ containing $q$ by Bertini's theorem \cite[Thm. 6.10]{jouanolou}, 
    so the geometric (or topological) genus is also 
    $p_{g}(F^{-1}(l))=\frac{1}{2}(d-1)(d-2)$ \cite[III, Remark 7.12.2]{AG}.

    There is an open subset $U \subset Y - \bigcup_{i} H_{i}$ such 
    that the restriction map $F|_{U}:U \to \ptwo$ is \'etale 
    \cite[I, Prop. 3.8]{Milne}.
    Removing a line $L \subset \ptwo$ gives an \'etale morphism
    $F: X = U - \pi_{p}^{-1}(L) \to \mathbb C^{2}$ of degree $d$ such that 
    the pre-image of every line meeting $F(X)$ is an irreducible curve.
    For the general such line $l$, the pre-image $F^{-1}(l)$ is a Zariski
    open subset of a smooth projective plane curve ($Y \cap H$ in the
    discussion above) of geometric genus
    $\frac{1}{2}(d-1)(d-2) > 0$. 
\end{example}

\begin{example}\label{topological}
    Let $S \subset \mathbb R^{3}$ be the cubic surface 
    $x^{2}+y^{2}=z^{2}(1-z)$ obtained by rotating a nodal cubic curve 
    about the $z$-axis and consider the projection $F:S \to \mathbb R^{2}$ 
    to the $xy$-plane. 
    The Jacobian determinant vanishes along the set $W$ consisting of the origin 
    and the circle $z=\frac{2}{3}$, so the restricted map 
    $X := S-W \to \mathbb R^{2}$ becomes a local diffeomorphism. 
    The map is not globally injective (points on the disk 
    $D=\{(x,y): x^{2}+y^{2} \leq \frac{4}{27}\}$ have three pre-images) and the 
    pre-image of a line $l$ is connected if and only if $l$ misses the 
    disk $D$ (compare \cite[Example 3.3]{NX}).
   
    Now we complexify the above example. Take $S \subset \mathbb C^{3}$ to be 
    the {\it complex} surface of the same equation and 
    remove the locus $W$ consisting of the origin and the 
    complex curve $z=\frac{2}{3}$ to obtain $X=S-W$. Projection to 
    the $xy$-plane now gives a non-injective local biholomorphism 
    $F: X  \to \mathbb C^{2}$. 
    Consider the complex line $l:y=ax$ through the origin: the 
    equations of the pre-image under $F$ are $y=ax$ and 
    $(1+a^{2})x^{2}=z^{2}(1-z)$. If $a \neq \pm i$, then the 
    pre-image $F^{-1}(l)$ is a connected rational curve 
    (a parametrization is given by $x=\frac{t^{2}-1-a^{2}}{t^{3}}, \; y= ax, \; z=tx$),
    but for $a = \pm i$ the pre-image is not connected, consisting of 
    two disjoint lines $z=0$ and $z=1$.
    The pre-images of the lines through the point $(0,1)$ are all 
    connected (being zero sets of irreducible polynomials), but 
    most of them are the complement of a smooth elliptic cubic curve 
    by finitely many points.  
\end{example}

\begin{example}\label{rigid}
Theorem \ref{main} fails if we remove the rigidity condition. 
Consider the map $F: D \times D \to \mathbb C^{2}$ given 
by $F(z,w)=((w+2) e^{7z}, w)$, where $D \subset \mathbb C$ is the 
open unit disc. 
Then $F$ is a local biholomorphism since $\det JF = 7(w+2) e^{7z} \neq 0$ on 
$D \times D$. 
Now let ${\mathcal {U}} \subset D \times D$ be the union of all 
connected components through $(0,0)$ of all pre-images of complex lines $l$ 
passing through $(2,0)$. 
It follows from the inverse function theorem that the interior 
${\mathcal U}^{0}$ of ${\mathcal U}$ is closed in ${\mathcal U}$. 
Since $(0,0) \in {\mathcal U}^{0}$ and ${\mathcal U}$ is connected, one has that 
${\mathcal U} = {\mathcal U}^{0}$ is open.  
The restriction $F|_{{\mathcal U}}$ is not injective, since 
$F(0,0)=F(\frac{2 \pi i}{7}, 0)=(2,0)$. 
The line $w=0$ pulls back to $D \times \{0\}$, a non-rigid domain in 
$\mathbb C \mathbb P^{1}$. 
Restricting $F$ to a sufficiently small neighborhood 
${\mathcal U}^{\prime} \subset {\mathcal U}$ of $D \times \{0\}$, 
it is clear that the pull-back of every line $l$ through $(2,0)$ is a 
connected set which is conformal to an open non-rigid subset of 
$\mathbb C \mathbb P^{1}$. 
\end{example}

%\pagebreak

\section{The continuity proof}

In this section we complete the proof of Theorem \ref{main} 
by showing continuity of the map 
$v:\mathbb C \mathbb P^{n-1} \to \mathbb R \mathbb P^{2n-1}$ there,  
constructed as follows. 
Starting with a local biholomorphism $F:X \to \mathbb C ^{n} $ such that 
$z_{1} \neq z_{2}$ and $F(z_1)=F(z_2)=0$, we assumed that for each complex 
line $l\in \mathbb C \mathbb P^{n-1}$ through the origin, $F^{-1}(l)$ is a rigid 
domain in $\mathbb C \mathbb P^{1}$, hence a connected rational curve.
Fixing $w_{1}=0$ and $w_{2}=1$ in $\mathbb C \mathbb P^{1}$ and real 
tangent vectors $\alpha_{i} \in T_{\mathbb P^{1},w_{i}}$, we showed that 
there are exactly two embeddings 
$T_l^1, T_{l}^{2}:F^{-1}(l) \to \mathbb C \mathbb P^1$ sending $z_i$ to 
$w_i$ ($i=1,2$) and satisfying 
$$dF(z_1)(dT_l^j(z_1))^{-1}\alpha_1=
dF(z_2)(dT_l^j(z_2))^{-1}\alpha_2 \in \mathbb C^n$$
for $j=1,2$.
Furthermore, both embeddings give rise to the same {\it real} line through $0$ 
(we now drop the superscript and write $T_l$ for either of the two embeddings).
In particular, we have defined a map 
$v:\mathbb C\mathbb P^{n-1} \to \mathbb R \mathbb P^{2n-1}$.
such that $v(l)=p(\kappa(l))$, where 
$p:\mathbb C^n - \{0\} \to \mathbb R \mathbb P^{2n-1}$ is the natural 
projection and $\kappa:\mathbb C\mathbb P^{n-1} \to \mathbb C^n -\{0\}$ is 
given by $\kappa(l)=dF(z_1)(dT_l(z_1))^{-1}\alpha_1$. 
We will now show that $v$ is continuous; in other words, for any 
convergent sequence $l_{q} \to l$ in $\mathbb C \mathbb P^{n-1}$, 
we will show that $v(l_{q}) \to v(l)$.

The salient feature is that the submanifolds $F^{-1}(l_{q})$
converge to $F^{-1}(l)$ uniformly over compact subsets of $X$.
Since we are proving that $v(l_{q}) \to v(l)$ for {\it any} sequence
$l_q \to l$, it suffices to show that $v(l_{q})$ has a {\it subsequence} 
converging to $v(l)$.
Indeed, this shows that $v(l)$ is the only accumulation point.
We will achieve this by essentially (though not literally) writing the map
$T_{l}:F^{-1}(l) \to \mathbb C \mathbb P^1$ as a limit of the maps $T_{q}=T_{l_{q}}$.
To this end, let $S_{q}$ be a sequence in the unitary group $U(n)$ such that
$S_{q} \rightarrow I_{n}$ and $S_{q}(l)=l_{q}$.

Fix a complete Riemannian metric on $X$. Since $F^{-1}(l)$ is connected,
there is a closed geodesic ball $K_{0}={\overline {B_{R}(z_{1})}}$ of some
radius
$R > 0$ centered at $z_{1}$ such that $z_{1},z_{2}$ lie in
the same connected component $C_{0}$ of the compact domain $K_{0} \cap
F^{-1}(l)$ for the Riemann surface $F^{-1}(l)$. Setting
$K_{m} = {\overline {B_{R+m}(z_{1})}}$, we obtain a countable increasing
exhaustion $\{K_{m}\}_{m \geq 1}$ of $X$. Clearly $z_{1},\;z_{2}$  lie in the
same connected component $C_{m}$ of $K_{m} \cap F^{-1}(l)$ for each $m \geq 1$.
Since $K_{m}$ is compact, there is $\delta=\delta(m) > 0$ such that
\begin{enumerate}
         \item $F|_{B_{\delta}(z)}$ is a biholomorphism onto its image for all
         $z \in K_{m}$.
         \item $F|_{B_{\delta}(z) \cup B_{\delta}(z^{\prime})}$ is injective
whenever $B_{\delta}(z) \cap B_{\delta}(z^{\prime})$ is nonempty.
\end{enumerate}
Indeed, for $x \in K_{m}$, there exists $\delta_x > 0$ such that
$F|_{B_{\delta_x}(x)}$ is a biholomorphism onto its image.
The open cover $\{B_{\frac{\delta_x}{4}}(x)\}_{x \in K_{m}}$
for $K_{m}$ has a finite subcover $\{B_{\frac{\delta_k}{4}}(x_k)\}_{k=1}^{r}$
and we set
$\displaystyle \delta=\delta(m)=\min_{1 \leq k \leq r}\{\frac{\delta_k}{4}\}$.
For $z \in B_{\delta}(x_{k}) \cap K_{m}$, we have
$B_{\delta}(z) \subset B_{2 \delta}(x_{k})$ which implies
condition 1 since $F$ takes the larger ball biholomorphically onto
its image. If $B_{\delta}(z^{\prime})$ meets $B_{\delta}(z)$,
then $B_{\delta}(z^{\prime}) \subset B_{4 \delta}(x_{k})$ so that
${B_{\delta}(z) \cup B_{\delta}(z^{\prime})} \subset B_{4 \delta}(x_{k})$
and the restriction of $F$ to the latter ball is injective, which
checks condition (2). Setting
$\displaystyle \delta_{m}=\min_{k \leq m} \delta(k)$ we obtain a
nonincreasing sequence with the same properties.

Now fix $m \geq 1$. The set 
$$S=\{(a,z) \in X \times (F^{-1}(l) \cap K_m):d(a,z)=\frac{\delta_m}{2}\}$$
is compact, hence the continuous function
$G(a,z)=|F(a)-F(z)|$ achieves its minimum value $\lambda$ on $S$; 
moreover $\lambda > 0$ because $F|_{B_{\delta_m}(z)}$ is
injective for each $z \in F^{-1}(l) \cap K_m$.
For $\mu = \frac{\lambda}{2}$ we have 
$B_{\mu}(F(z)) \subset F(B_{\delta_{m}}(z))$ 
for all $z \in F^{-1}(l) \cap K_{m}$.
The convergence $S_{q} \to I$ yields $N_{m}>0$ such that the absolute values of 
the entries of the matrix $S_{q}-I_{n}$ are all less than 
$\frac{\mu}{n^{2} M}$ for $q \geq N_{m}$, where
$$ M=M(m) = \sup_{z \in F^{-1}(l) \cap K_{m}} |F(z)|.$$
With these choices, we have that $|S_{q}(F(z))-F(z)| < \mu$ for 
$z \in F^{-1}(l) \cap K_{m}$ and $q \geq N_{m}$ so that 
$S_{q}(F(z)) \in F(B_{\delta_{m}}(z))$.
As above, we may assume that the $N_{m}$ are nondecreasing.

For $\delta_{m}$ and $N_{m}$ as above, we have thus constructed 
well-defined holomorphic injections 
$\phi_{m,q} : C_{m} \hookrightarrow F^{-1}(l_{q})$ for $q \geq N_{m}$ given by
$$\phi_{m,q}(z)=[(F|_{B_{\delta_{m}(z)}})^{-1} \circ S_q \circ F](z).$$
By construction it is clear that $\phi_{m,q}$ fixes $z_1, z_2$
and that $\phi_{m^{\prime},q}|_{C_{m}}=\phi_{m,q}$ for
$m^{\prime} \geq m$ and $q \geq N_{m}$. 
Injectivity follows from condition (2) on $\delta(m)$ above. 
Composition with $T_{q}$ yields injective holomorphic maps
$$\psi_{m,q}=T_q \circ \phi_{m,q}: 
C_{m} \hookrightarrow \mathbb C \mathbb P^{1}$$ 
for $q \geq N_{m}$ which satisfy $\psi_{m,q}(z_{i})=w_{i}$ and 
$\psi_{m^{\prime},q}|_{C_{m}}=\psi_{m,q}$ for $m^{\prime} \geq m$. 
Note that we will make a convenient choice of identification 
$\mathbb C \mathbb P^{1} = \mathbb C \cup \{\infty\}$ later in the proof.

We will need the following compactness result in the proof of Lemma
\ref{subseq} below:

\begin{lemma}\label{normal} Let $D \subset \mathbb C$ be the open unit disc,
$0\neq a \in D$ and $0\neq b \in \mathbb C$.
Then 
$$
\mathcal F_{a,b}=\{\mbox{holomorphic injections } 
f:D \hookrightarrow \mathbb C 
\;|\;f(0)=0 \mbox{ and } f(a)=b\}
$$
is a normal family of functions.
\end{lemma}

\begin{proof} Any univalent holomorphic function $f$ on $D$ that
vanishes at $0$  satisfies the well-known estimate of Bieberbach
and Gronwall (see, for example \cite[Theorem 1.3]{hayman} for the
normalized estimate)
$$
\frac{|f'(0)|r}{(1+r)^2} \leq |f(z)| \leq
\frac{|f'(0)|r}{(1-r)^2},
$$
where $r=|z|$. Evaluating the inequality on the left at $z=a$ gives
an upper bound on $|f^{\prime}(0)|$. The inequality on the right now shows
that
$|f(z)|$ is uniformly bounded over compact subsets of $D$, with a
bound independent of $f$. By Montel's theorem $\mathcal F_{a,b}$ is normal,
that is ${\overline {\mathcal F_{a,b}}}$ is compact in the space of holomorphic
functions on $D$, endowed with the topology of uniform convergence over
compact subsets.
\end{proof}
\vskip5pt

\begin{lemma}\label{subseq}
         Fix $m \geq 1$ and an identification 
	 $\mathbb C \mathbb P^{1}=\mathbb C \cup \{\infty\}$. 
	 Let $\{\varphi_{k}\}$ be a subsequence of 
	 $\{\psi_{m,q}\}_{q \geq N_{m}}$ such that each $\varphi_{k}$ has 
	 a pole $p_{k}$ and $\{p_{k}\}$ converges to a point $x$ in the interior 
	 of $C_{m}$. 
	 Then $\{\varphi_{k}\}$ has a further subsequence converging 
	 uniformly over compact subsets of the interior of $C_m-\{x\}$ to a 
	 holomorphic injection of the interior of $C_m-\{x\}$ into 
	 $\mathbb C=\mathbb C \mathbb P^1-\{\infty\}$.
\end{lemma}

\begin{proof}
To set notation, write $\varphi_{k}=\psi_{m,s_{0}(k)}$ for some strictly 
increasing function $s_{0}:\mathbb N \to \mathbb N$ with $s_{0}(1) \geq N_{m}$.
Cover $C_{m} - \{x\}$ with countably many embedded holomorphic 
disks $\xi_{r}: D \hookrightarrow C_{m} - \{x\}$ in such a way that
$\xi_{r}(0) = z_{1}, \; z_{2} \in \xi_{r}(D)$ but 
$x \not\in \overline{\xi_{r}(D)}$.

Consider the $r$th embedded holomorphic disk 
$\xi_{r}: D \hookrightarrow C_{m} - \{x\}$.
Then $p_{k} \not\in \xi_{r}(D)$ for $k$ sufficiently large because 
$x \not\in \overline{\xi_{r}(D)}$, hence the maps 
$\varphi_{k}\circ\xi_{r}$ belong to class 
$\mathcal F _{\xi_{r}^{-1}(z_2), w_2}$ of Lemma \ref{normal} (their 
images lie in $\mathbb C$ via the identification 
$\mathbb C \mathbb P^{1}=\mathbb C \cup \{\infty\}$). 
In particular, there exists a subsequence 
$\varphi_{s_{1}(k)}\circ\xi_{r}$ 
converging uniformly over the compact sets of $D$ by normality of 
$\mathcal F _{\xi_{r}^{-1}(z_2), w_2}$.
Composing with $\xi_{r}^{-1}$, the subsequence $\varphi_{s_{1}(k)}$ converges 
uniformly over the compact subsets of the topological disc
$\xi_{r}(D) \subset C_{m}- \{x\}$ to a holomorphic map taking
values in $\mathbb C = \mathbb C \mathbb P^{1} - \{\infty\}$.
Since $\varphi_{s_{1}(k)} (z_i)=w_i \;\;(i=1,2)$, we see that the nonconstant
limit function is the local uniform limit of a sequence of injective
holomorphic functions. 
By a well-known theorem of Hurwitz, the limit function is itself injective.

Taking $r=1$ above, $\varphi_{k}$ has a
subsequence $\varphi_{s_{1}(k)}$ which converges uniformly over compact
sets in $\xi_{1}(D)$ ($s_{1}:\mathbb N \to \mathbb N$ strictly
increasing as above). Applying the previous paragraph with $r=2$ and
the subsequence $\varphi_{s_{1}(k)}$, we obtain a further subsequence
$\varphi_{s_{2}(k)}$ with the same property with respect to both 
$\xi_{1}(D)$ and $\xi_{2}(D)$. Continuing in this fashion,
we obtain a countable sequence of subsequences $\varphi_{s_{r}(k)}$
(each uniformly convergent on compact subsets of $\xi_{p}(D)$ for $1 
\leq p \leq r$), each being extracted from the previous:
$s_{1}(\mathbb N) \supset  \; s_{2}(\mathbb N) \supset \; 
s_{3}(\mathbb N) \dots$.
Thus the diagonal subsequence $\varphi_{s_{r}(r)}$ converges
uniformly to an injective holomorphic function on every compact subset of
$C_{m}-\{x\}$. This proves Lemma \ref{subseq}.
\end{proof}

\vskip5pt

Taking $m=1$, we now make a choice of both $\infty \in \mathbb C 
\mathbb P^{1}$ and a subsequence $\varphi_{k}$ of $\psi_{1,q}$ for which 
the hypotheses of Lemma \ref{subseq} holds. 
Start with any point $a \in \mathbb C \mathbb P^{1}$. 
If there happens to be a subsequence $\psi_{1,q_{k}}$ and $p_{k} \in C_{1}$ with 
$\psi_{1,q_{k}}(p_{k})=a$ and $p_{k} \to x$ with $x$ in the interior of 
$C_{1}$, then we simply take $\infty=a$ and $\varphi_{k}=\psi_{1,q_{k}}$. 
The alternatives are that 
(1) $a \in \psi_{1,q}(C_{1})$ for only {\it finitely} many $q$ or 
(2) $a \in \psi_{1,q}(C_{1})$ for infinitely many $q$, but for all 
subsequences $p_{k}, \; q_{k}$ with $\psi_{1,q_{k}}(p_{k})=a$ and $p_{k} \to x$, 
we have $x \in \partial C_{1}$. 
In either case, we can find a subsequence $q_{k}$ and an 
embedded holomorphic disk $\xi: D \hookrightarrow C_{1}$ such that 
$\xi(0)=z_{1}$, $z_{2} \in \xi(D)$ and $a \not\in \psi_{1,q_{k}}(\xi(D))$ 
for $k$ sufficiently large. 
Thus 
$\psi_{1,q_{k}} \circ \xi :D \hookrightarrow \mathbb C 
= \mathbb C \mathbb P^{1} - \{a\}$ 
for $k > \; > 0$ and by Lemma \ref{normal} we can find a further subsequence 
$\psi_{1,q_{k_{l}}} \circ \xi$ converging uniformly over compact subsets. 
Now fix $y \in D$, set $x = \xi(y)$ and take 
$\infty=\lim_{l \to \infty} \psi_{1,q_{k_{l}}}(x)$. 
Thus we obtain a subsequence $\varphi_{k}=\psi_{1,q_{k_{l}}}$ of 
$\psi_{1,q}$ for which every further subsequence forcibly has poles which 
converge to $x$ in the interior of $C_{1}$ and we may apply Lemma 
\ref{subseq}. 

Finally we vary $m \geq 1$. 
Starting with $m=1$, choose $\infty \in \mathbb C \mathbb P^{1}$ and 
$\varphi_{k}$ as above and apply Lemma \ref{subseq}
to obtain a subsequence $\varphi_{s_{1}(k)}$ which converges uniformly 
to an injective holomorphic function on compact subsets of the interior of 
$C_{1}-\{x\}$. 
Setting $B = \min\{k:s_{1}(k) \geq N_{2}-1\}$ and $t_{1}(k)=s_{1}(k+B)$ 
we have a subsequence $\varphi_{t_{1}(k)}$ of $\psi_{2,n}$ whose
restrictions to $C_{1}$ form a subsequence of $\varphi_{s_{1}(k)}$; 
moveover, the hypotheses of Lemma \ref{subseq} hold automatically by the 
compatibility condition on the maps $\psi_{m,q}$, that is, 
$\psi_{m^{\prime},q}|_{C_{m}}=\psi_{m,q}$ for $m^{\prime} \geq m$. 
Thus applying Lemma \ref{subseq} with $m=2$ gives a further subsequence that
converges uniformly to an injective holomorphic function on compact
subsets of the interior of $C_{2}-\{x\}$. 
Continuing in this fashion, we obtain a countable sequence 
$\{\varphi_{s_{m}(k)}\}_{m \geq 1}$ of subsequences, each being extracted 
from the restriction of the next: 
$\varphi_{s_{m}(k)}$ converges uniformly over compact subsets to an 
injective holomorphic function on the interior of $C_m-\{x\}$.
Taking the diagonal subsequence $\varphi_{s_{m}(m)}$ yields

\begin{lemma}\label{final}
         There is an injective holomorphic map
         $\eta: F^{-1}(l)-\{x\} \hookrightarrow \mathbb C=
	 \mathbb C \mathbb P^1-\{\infty\}$ 
	 such that for each fixed $m \geq 1$, the restriction of $\eta$ to the
         interior of $C_m-\{x\}$ is the local uniform limit of a 
	 subsequence $\psi_{m,q_{l}}$ of $\psi_{m,q}$.
\end{lemma}

The map $\eta$ extends to a holomorphic injection 
$\tilde{\eta}:F^{-1}(l)\to \mathbb C \mathbb P^1$ such that
$\tilde{\eta}(z_i)=w_i$ for $i=1,2$. 
For $z$ in the interior of $\in C_1-\{x\}$ 
(say $z$ near $z_1$ or $z_2$), we use Lemma \ref{final} to write 
$$\tilde{\eta}(z)=\lim_{k\to\infty} \psi_{1,q_{k}}(z)$$
as a local uniform limit for some subsequence $q_{k}$.
Since the sequence of derivatives also converges, we have
$$dF(z_1) \circ (d\tilde{\eta}(z_1))^{-1}(\alpha_1) =\lim_{k\to \infty}
dF(z_1) \circ (d\psi_{1, q_k}(z_1))^{-1}(\alpha_1).$$
Since $\psi_{1,q}=T_q\circ \phi_{1,q}$ and
$\phi_{1,q}(z)=[(F|_{(B_{\delta_{1}(z)})})^{-1}\circ S_q \circ F](z)$, 
this becomes
$$ \lim _{k\to \infty} [ (S_{q_k})^{-1} \circ dF(z_1) \circ
(dT_{q_k}(z_1))^{-1}(\alpha_1)],$$
but $\lim_{k \to \infty} S_{q_{k}} = I=n \times n$ identity 
matrix, so we conclude that
$$ dF(z_1) \circ (d\tilde{\eta}(z_1))^{-1}(\alpha_1)
=\lim_{k\to \infty} dF(z_1) \circ (dT_{q_k}(z_1))^{-1}(\alpha_1). $$
Analogously, $$ dF(z_2) \circ (d\tilde{\eta}(z_2))^{-1}(\alpha_2)=\lim_{k\to \infty} dF(z_2) \circ
(dT_{q_k}(z_2))^{-1}(\alpha_2). $$
By the definition of $T_{q_k}$, one has 
$$dF(z_1) \circ (dT_{q_k}(z_1))^{-1}(\alpha_1)=
dF(z_2) \circ (dT_{q_k}(z_2))^{-1}(\alpha_2),$$
so that 
$$ dF(z_1) \circ (d\tilde{\eta}(z_1))^{-1}(\alpha_1)=
dF(z_2) \circ (d\tilde{\eta}(z_2))^{-1}(\alpha_2).$$
Together with $\tilde{\eta}(z_i)=w_i$ ($i=1,2$), the last relation shows that 
$\tilde{\eta}$ is equal to one of the two maps $T_{l}^{j}$. 
Taking images in $\mathbb R \mathbb P^{2n-1}$ one now sees that 
$v(l_{q_k}) \to v(l)$. 
This shows that the section 
$v:\mathbb C \mathbb P^{n-1} \to \mathbb R \mathbb P^{2n-1}$ of the natural map 
$\pi: \mathbb R \mathbb P^{2n-1} \to \mathbb C \mathbb P^{n-1}$ is continuous. 
Following the earlier part of the proof, this contradiction proves 
Theorem \ref{main}.
\qed

\end{document}